%This file starts with my \define commands,
%most of which are probably not used 
%in any given file.
%Hopefully they do not cause any trouble.

%The file uses AMSTEX and AMSppt.
%Depending on how your system is set up,
%you may need one of the following commands:

\input amstex
\documentstyle{amsppt}

\magnification\magstep1

\define\p{\Bbb P}

\define\a{\Bbb A}

\redefine\c{\Bbb C}

\redefine\o{\Cal O}

\define\q{\Bbb Q}

\define\r{\Bbb R}

\define\z{\Bbb Z}

\define\map{\dasharrow}

\define\qtq#1{\quad\text{#1}\quad}

\define\section#1{

\bigpagebreak{\smc#1}\bigpagebreak

}

\define\rup#1{\ulcorner #1\urcorner}

\define\rdown#1{\llcorner #1\lrcorner}

\define\demop{\demo{Proof}}

\define\pic{\operatorname{Pic}}

\define\spec{\operatorname{Spec}}

\define\red{\operatorname{red}}

\define\len{\operatorname{length}}

\define\im{\operatorname{im}}

\define\chr{\operatorname{char}}

\define\deq{:=}

\define\norm#1{\vert\vert#1\vert\vert}

\define\broot#1{[@,@,@,\root{n}\of{#1}@,@,@,]}

\define\bir{\operatorname{Bir}}

\topmatter
\title     Nonrational Covers of  $\c\p^m\times \c\p^n$
\endtitle
\author J\'anos Koll\'ar
\endauthor
 \endtopmatter

\head 1. Introduction
\endhead

The aim of this article is to provide further examples of higher
dimensional varieties which are rationally connected, but not
rational and not even ruled. The original methods of 
\cite{Iskovskikh-Manin71; Clemens-Griffiths72}  were further
developed    by many authors (for instance, 
\cite{Beauville77; Iskovskikh80; Bardelli84}) and they give a quite
complete picture in dimension three.  \cite{Corti96} is a recent
overview. For a while only special examples have been known in higher
dimensions 
\cite{Artin-Mumford72; Sarkisov81,82; Pukhlikov87; CTO89}.  For
hypersurfaces in $\p^n$ the rationality question was considered in
\cite{Koll\'ar95}. There I proved the following:

\proclaim{1.1 Theorem} \cite{Koll\'ar95}  Let $X_d\subset \c\p^{n+1}$
be a very general hypersurface of degree $d$. Assume that
$$
\frac{2}{3}n+3\leq d\leq n+1.\tag 1.1.1
$$ Then $X_d$ is not rational.\qed
\endproclaim

 Here ``very general" means that the result holds for hypersurfaces
corresponding to a point in the complement of countably many closed
subvarieties in the space of all hypersurfaces.

The method can be applied to hypersurfaces in $\p^m\times \p^{n+1}$.
Let
$X_{c,d}$ be such a hypersurface of bidegree $(c,d)$. Via the
projection
$X_{c,d}\to \p^m$ it can be viewed as  a family of degree $d$
hypersurfaces in $\p^{n+1}$ parametrized by $\p^m$.   $X_{c,d}$ is
rationally connected if
$d\leq n+1$, no matter what $c$ is (cf. \cite{Koll\'ar96, IV.6.5}).

A straightforward application of the method of \cite{Koll\'ar95} 
provides an analog of (1.1) for hypersurfaces when  $c\geq m+3$ and
$d,n$ satisfy the inequalities   (1.1.1). It is, however, of more
interest to study some cases where the fibers of
$X_{c,d}\to \p^m$ are rational. Here I propose to work out two cases:
conic bundles and families of cubic surfaces. 

The method of \cite{Koll\'ar95} works naturally for cyclic covers,
and it is easier to formulate the results that way.

Fix a prime $p$ and let  $X_{ap,bp}\to \c\p^m\times \c\p^n$ be a
degree
$p$ cyclic  cover ramified along  a very general hypersurface of
bidegree
$(ap,bp)$.   \cite{Koll\'ar95} shows that $X_{ap,bp}$ 
 is not rational   and not even ruled if $ap>m+1$ and $bp>n+1$.  In
this paper I study the case when $bp=n+1$ and $n=1,2$.  The main
results are the following:

\proclaim{1.2 Theorem}   Let $X_{2a,2}\to \c\p^m\times \c\p^1$ be a
double cover ramified along  a very general hypersurface of bidegree
$(2a,2)$; $m\geq 2$. Then
$X_{2a,2}$ is not rational  if $2a> m+1$.

More precisely, if $Y$ is any variety of dimension $m$ and
$\phi:Y\times
\p^1\map X_{2a,2}$  a dominant map then $2|\deg \phi$. 
\endproclaim

\demo{1.2.1 Remarks}  (1.2.1.1) $X_{2a,2}\to \c\p^n$ is a conic
bundle. For conic bundles we have the very strong results of 
\cite{Sarkisov81,82}  which say that a conic bundle is not rational
if the locus of singular fibers plus 4 times the canonical class of
the base is effective. In our case the locus of singular fibers is a
divisor in $\p^m$ of degree $4a>2m+2$.  This is  lower than the
Sarkisov bound  $4m+4$. Sarkisov's normal crossing assumptions are
also not satisfied by
$X_{2a,2}\to \c\p^n$.  Thus some of our cases are not covered by the
results of \cite{Sarkisov81,82}. This suggests the possibility that 
the bounds of \cite{Sarkisov81,82} can be improved considerably. This
is not clear even for conic bundles over surfaces.  

(1.2.1.2) It is natural to ask if the projection $X_{2a,2}\to
\c\p^m$ is the only conic bundle structure of $X_{2a,2}$ or not. The
proof gives such examples in characteristic 2, but not over $\c$.
\enddemo

For families of cubic  surfaces the result is   weaker:

\proclaim{1.3 Theorem}  Let $X_{3a,3}\to \c\p^m\times \c\p^2$ be a
cyclic triple cover ramified along  a very general hypersurface of
bidegree $(3a,3)$;
$m\geq 1$. Then $X_{3a,3}$ is not rational   and not even ruled if
$3a> m+1$.
\endproclaim

\demo{1.3.1 Remark}    Similar results for $m=1$ were proved by
\cite{Bardelli84}. 
\enddemo

For hypersurfaces in $\c\p^m\times \c\p^{n+1}$ these imply the
following:

\proclaim{1.4 Theorem}  (1.4.1)  Let $X_{c,2}\subset \c\p^m\times
\c\p^2$ be a very general hypersurface of bidegree $(c,2)$; $m\geq
2$. Then $X_{c,2}$ is not rational  if $c\geq m+3$.

(1.4.2)  Let $X_{c,3}\subset \c\p^m\times \c\p^3$ be a very general
hypersurface of bidegree $(c,3)$; $m\geq 1$. Then $X_{c,3}$ is not
rational     if $c\geq m+4$.
\endproclaim

\demop  By \cite{Koll\'ar96, V.5.12--13} it is sufficient to find an
algebraically closed  field
$k$ and a single example of a nonruled hypersurface
$X_{c',2}\subset \p^m\times \p^2$ resp. 
$X_{c',3}\subset \p^m\times \p^3$ over $k$ for some $c'\leq c$. We
proceed to construct such  examples in characteristic two for conic
bundles and in characteristic three for families of cubic surfaces.

Consider $\a^m\times\a^n$ with coordinates $(u_1,\dots,u_m;
v_1,\dots,v_n)$.  Let $f(u,v)$ be a polynomial of bidegree $(c,n+1)$. 
 We have a hypersurface
$$ Z\deq (y^{n+1}-f(u,v)=0)\subset V\deq \a^1\times \a^m\times\a^n,
$$ where $y$ is a coordinate on $\a^1$. There are several natural
ways to associate  a projective variety to $Z$:

(1.4.3.1) We can view  $V$ as a coordinate chart in $\p^m\times
\p^{n+1}$. The closure $\bar Z_1$ of $Z$ is  a hypersurface of 
 bidegree  $(c,n+1)$.  For $n=1,2$ this gives our examples $X_{c,2}$
and
$X_{c,3}$.

(1.4.3.2)  We can compactify $\a^m\times\a^n$ to
$\p^m\times\p^n$ and view $f$ as a section of the line bundle 
$\o_{\p^m\times\p^n}(c,n+1)$.  Assume that $n+1|c$ and let
$L=\o(c/(n+1),1)$.  The corresponding cyclic cover  $\bar
Z_2=\p^m\times\p^n[\root{n+1}\of{f}]$  (defined in (2.3)) gives
another compactification of
$Z$. 

(1.2--3) show using the second representation that  $\bar Z_2$ is not
ruled if
$f$ is very general and $c>n+1$. This implies (1.4).\qed\enddemo

\demo{1.5 Generalizations}  It is clear from the proof that it
applies to many different cases. The main point is to have a family
of conics or cubic surfaces   whose branch divisor is sufficiently
large, but I found it hard to write down a reasonably general
statement.
\enddemo

\demo{1.6 Outline of the proof of (1.2--3)}  By \cite{Koll\'ar96,
V.5.12--13} it is sufficient to find an algebraically closed  field
$k$ and a single example of a nonruled cyclic cover 
$X_{2a',2}\to \p^m\times \p^1$   (resp. $X_{3a',3}\to \p^m\times
\p^2$ ) for some $a'\leq a$.  We proceed to construct such  examples
in characteristic two for conic bundles  and in characteristic three
for families of cubic surfaces.

Fix an algebraically closed field $k$ of characteristic $p$. Let
$\pi: X_{ap,bp}\to \p^m\times \p^n$ be a degree
$p$ cyclic  cover ramified along  a very general hypersurface of
bidegree
$(ap,bp)$.  
$\pi$ is purely  inseparable, thus $X_{ap,bp}$ is (purely 
inseparably) unirational. We intend to show that it is frequently not
ruled. 

Section 2 is a    review of the  machinery of inseparable cyclic
covers and their applications to nonrationality problems.
  $X_{ap,bp}$ has isolated singularities by (2.1.4) and (2.2.3). 
These can be resolved, and if
$\pi_Y: Y\to X_{ap,bp}\to \p^m\times \p^n$ is a resolution then by
(2.5) there is a nonzero map
$$ \pi_Y^*\o(ap-m-1,bp-n-1)\to
\wedge^{m+n-1}\Omega_Y^1\cong\Omega_Y^{m+n-1}.
$$ Thus  \cite{Koll\'ar95} shows that $X_{ap,bp}$ 
 is not rational   and not even ruled if $ap>m+1$ and $bp>n+1$.  In
this paper I study the case when $ap>m+1$ and $bp=n+1$.

Section 3 contains a nonruledness criterion.  Let $f:Y\to \p^m$ be
the composite of $\pi_Y$ with the first projection.  (3.2) shows that
if $ap>m+1$ and $bp=n+1$, then 
$Y$ is ruled iff the generic fiber of $f$ is ruled over the function
field
$k(\p^m)$.  This is the main technical departure from
\cite{Koll\'ar95}. There I used varieties in positive characteristic
which were shown to be not even separably uniruled. Here the
varieties in question are separably uniruled, but we are able to get
a   description of all separable unirulings. 

Rationality or ruledness over nonclosed fields is a very interesting
question. Unfortunately I can not say much, except for $n=1,2$. 

When $n=1$, the generic fiber of $f$ is a plane conic. Conics over
nonclosed fields are considered in section 4. We get a complete
description of their unirulings in  (4.1). Theorem (1.2) is implied
by (2.5), (3.1) and (4.2). 

If $n=2$, the generic fiber of $f$ is a cubic surface $S$ given by an
equation
$u^3=f_3(x,y,z)$. Since $X_{ap,bp}$ has only  isolated singularities,
$S$ is nonsingular over $k(\p^m)$, though we will see that it is not
smooth (5.2).  Thus we are led to investigate  nonsingular Del Pezzo
surfaces over arbitrary fields  in section 5. We are forced to study
the situation over nonperfect fields; this introduces several new
features.  The main result is (5.7) which generalizes   results of
Segre and Manin to nonperfect fields.
\enddemo

\demo{1.7 Terminology} I follow the terminology of 
\cite{Koll\'ar96}. 

If a scheme $X$ is defined over a field $F$, and
$E$ is a field extension then $X_E$ denotes the scheme obtained by
base extension. 

For a field $F$, $\bar F$ denotes an algebraic closure. 

$X_F$ is called rational, if it is birational to $\p^n_F$, and the
birational map is defined over $F$. Sometimes for emphasis I say that
$X_F$ is rational over $F$. The same convention applies to other
notions (ruled, uniruled, irreducible etc.).

If $X_{\bar F}$ is rational then I say that $X_F$ is geometricaly
rational. Similarly for  other notions (ruled, uniruled, irreducible
etc.).

Following standard terminology, a scheme is called nonsingular if all
of its local rings are regular. Over nonperfect fields this is not
the same as being smooth. 
\enddemo

\demo{Acknowledgements}  Partial financial support was provided by the
NSF under grant number  DMS-9102866. These notes were typeset by
\AmSTeX, the \TeX macro system of the American Mathematical Society.
\enddemo

\head 2. Inseparable cyclic covers
\endhead

First we recall the definitions and basic properties of critical
points of sections of line bundles in positive characteristic. For
proofs see
\cite{Koll\'ar95} or \cite{Koll\'ar96, V.5}. 

\demo{2.1 Definition}  Let $X$ be a smooth variety over an
algebraically closed field
$k$ and
$f$ a function on
$X$. Let $x\in X$ be a closed point and assume that $f$ has a critical
point at $x$. Choose local coordinates
$x_1,\dots, x_n$ at $x$.  

(2.1.1)  $f$ has a {\it nondegenerate critical point} at $x$  iff
$\partial f/\partial x_1,
\dots, \partial f/\partial x_n$ generate the maximal ideal of the
local ring
$\o_{x,X}$.  This notion is independent of the local coordinates
chosen.

If $\chr k\neq 2$  or $\chr k=2$ and $n$ is even then $f$ has a 
nondegenerate critical point at $x$  iff  in suitable local
coordinates $f$ can be written as
$$ f=c+
\cases x_1x_2+x_3x_4+\dots+x_{n-1}x_n+f_3,\qtq{if $n$ is even,}\\
x_1^2+x_2x_3+\dots+x_{n-1}x_n+f_3,\qtq{if $n$ is odd,}\\
\endcases
\qtq{where} f_3\in m_x^3.
$$

(2.1.2) If $\chr k=2$ and $\dim X$ is odd, then every critical point
is degenerate. 

(2.1.3)  Assume that  $\chr k=2$ and  $\dim X$ is odd. A critical
point of $f$ is called {\it almost nondegenerate} iff
$\len \o_{x,X}/(\partial f/\partial x_1,\dots, \partial f/\partial
x_n)=2$. Equivalently, in suitable local coordinates $f$ can be
written as
$$ f=c+ax_1^2+x_2x_3+\dots+x_{n-1}x_n+bx_1^3+f_3\qtq{where $b\neq 0$.}
$$

(2.1.4)  Assume that $\chr k|d$. Then the hypersurface
$Z=(y^d-f(x_1,\dots, x_n)=0)$ is singular at the point $(y,x)\in Z$
iff
$x\in X$ is a critical point of $f$.

(2.1.5) Let $L$ be a line bundle on $X$ and  $s\in H^0(X,L^d)$  a
section. Let $U\subset X$ be an open affine subset  such that
$L|U\cong \o_U$. Choose such an isomorphism. Then $s|U$ can be viewed
as section of
$\o_U^{\otimes d}\cong
\o_U$. Thus it makes sense to talk about its critical points. If
$\chr k|d$ then this is independent of the choice of  $U$ and of the
trivialization
$L|U\cong \o_U$. (This  fails if the characteristic does not divide
$d$.)
\enddemo

The usual  Morse lemma can be generalized to positive characteristic.
We use
 it in a somewhat technical form. 

\proclaim{2.2 Proposition}  Let $X$ be a smooth variety over a field
of $\chr p$ and
$L$  a line bundle on $X$. Let $d$ be an integer divisible by $p$ 
and $W\subset H^0(X,L^d)$ a finite dimensional subvectorspace.  Let
$m_x$ denote the ideal sheaf of $x\in X$. Assume that: 

(2.2.1)  For every closed point
$x\in X$ the restriction map 
$W\to  (\o_X/m_x^2)\otimes L^k$ is surjective,

(2.2.2)  For every closed point
$x\in X$ there is an $f_x\in W$  which has an (almost) nondegenerate
critical point at $x$.

Then a general section $f\in W$ has only (almost) nondegenerate
critical points.
\endproclaim

\demop  This is a simple constant count.  Fix $x\in X$ and let
$W_x\subset W$ be the set of functions with a critical point at $x$.
By (2.2.1), $W_x$ has codimension $n$. In $W_x$ the functions with an
(almost) nondegenerate critical point at $x$ form an open subset
$W_x^0$ which is nonempty by (2.2.2).  Thus the set of functions with
a degenerate critical point is $\cup_x (W_x-W_x^0)$ and it has
codimension at least one in $W$. \qed
\enddemo

\proclaim{2.2.3 Lemma}  Let $X_1,X_2$ be smooth varieties  over a
field of $\chr p$ and $L_i$  very ample line bundles on $X_i$.  Let
$L=p_1^*L_1\otimes p_2^*L_2$ be the corresponding line bundle on
$X=X_1\times X_2$.  If $p|d$  then a general section
$f\in H^0(X,L)$ has only (almost) nondegenerate critical points.
\endproclaim

\demop Pick a point $x=(x_1,x_2)$.  The condition (2.2.1) is clearly
satisfied. In order to check (2.2.2) choose global sections $u_i\in
H^0(X_1,L_1)$ 
 and $v_j\in H^0(X_2,L_2)$ such that they give local coordinates at
$x_1$ resp. $x_2$. 

If $p\neq 2$ then $\sum u_i^2+\sum v_j^2$ gives a section of $L^d$
with a nondegenerate critical point at $x$.

If $p=2$ we need to consider a few cases. Set $n_i=\dim X_i$. We plan
to use the function 
$$ g=\sum_{1\leq i\leq n_i/2}u_{2i-1}u_{2i} + \sum_{1\leq j\leq n_2/2}
v_{2j-1}v_{2j}.
$$ In both of the $n_i$ are even, then we can take $f=g$.  If both of
the $n_i$ are odd, then we can use $f=g+u_{n_1}v_{n_2}$.  Otherwise
we may assume that
$n_1$ is odd and $n_2$ is even.  Then we use
$$ f=g+u_{n_1}v_{n_2-1}+u_{n_1}^2v_{n_2}.
$$ Explicit computation shows that $f$ has an almost nondegenerate
critical point.\qed\enddemo

\demo{2.3 Definition}  Cyclic covers

Let $X$ be a   scheme, $L$ a line bundle on
$X$ and $s\in H^0(X, L^d)$ a section.  Assume for simplicity that the
divisor of its zeros $(s=0)$ is reduced.  The cyclic cover of $X$
obtained by taking a $d^{th}$-root of $s$, denoted by
$X[\root{d}\of{s}]$ is a scheme locally constructed as follows:

Let $U\subset X$ be an open set such that $L|U\cong \o_U$. Then $s|U$
can be identified with a function $f\in H^0(U,\o_U)$. Let
$V\subset \a^1\times U$ be the closed subset defined by the equation
$y^d-f=0$ where $y$ is the coordinate on $\a^1$. The resulting
schemes can be patched together in a natural way to get a scheme
$X[\root{d}\of{s}]$; 
 cf. \cite{Koll\'ar96, II.6.1}. We are interested in it only up to
birational equivalence, so the precise definitions are unimportant.
\enddemo

The only result about cyclic covers we need is the following special
case of
\cite{Koll\'ar96, V.5.10}:

\proclaim{2.4 Proposition}  Let $X$ be a  smooth variety of dimension
$n$ over a field
$k$ of
$\chr p$,
$L$ a line bundle on
$X$ and $d$ an integer divisible by $p$. Let $s\in H^0(X, L^d)$ be a
section with  (almost) nondegenerate critical points.  Let $\pi: Y\to
X$  be a smooth projective model of
$X[\root{d}\of{s}]$ ($Y$ always exists). 

Then there is a nonzero map
$$ 
\pi^*(K_X\otimes L^d)\to
\wedge^{n-1}\Omega_Y^1\cong\Omega_Y^{n-1}.\qed
$$
\endproclaim

Applied to the cyclic covers $\bar
Z_2=\p^m\times\p^n[\root{p}\of{s}]$ from the introduction, we get the
following:

\proclaim{2.5 Corollary}   Fix a prime $p$ and let $k$ be an
algebraically closed field
 of characteristic 
$p$. Let $s\in H^0(\p^m\times\p^n,\o(a,b)^{\otimes p})$ be
 a general section and $q: Y\to \p^m\times\p^n$    a smooth
projective model of
$\p^m\times\p^n[\root{p}\of{s}]$. 

Then there is a nonzero map
$$ q^*\o(ap-m-1,bp-n-1)\to
\wedge^{m+n-1}\Omega_Y^1\cong\Omega_Y^{m+n-1}.\qed
$$
\endproclaim

\head 3. A nonruledness criterion
\endhead

In this section we prove the following generalization of
\cite{Koll\'ar96, V.5.11}.

\proclaim{3.1 Theorem} Let $X,Y$ be smooth proper varieties and
$f:Y\to X$ a surjective morphism, $n=\dim Y$. 
 Let $M$ be a big line bundle on $X$ and assume that for some $i>0$ 
there is a nonzero map
$$ h:f^*M\to \wedge^i\Omega_Y^1.
$$

(3.1.1)  Let $Z$ be an affine variety of dimension $n-1$ and
$\phi:Z\times
\p^1\to Y$ a dominant and separable morphism. Then there is a morphism
$\psi:Z\to X$ which fits in the commutative diagram
$$
\CD Z\times \p^1 @>\phi>> Y\\ @VVV   @VVfV\\ Z @>\psi>> X.
\endCD
$$

(3.1.2) Let $F=k(X)$ be the field of rational functions on $X$ and 
$Y_F$  the generic fiber of $f$.   There is a one-to-one
correspondence 
$$
\left\{\foldedtext\foldedwidth{2.8cm}{degree $d$ separable 
unirulings of
$Y$}\right\}\leftrightarrow
\left\{\foldedtext\foldedwidth{2.8cm}{degree $d$ separable unirulings
of $Y_F$ }\right\}.
$$ In particular, 
$Y$ is ruled  iff $Y_F$ is  ruled
 over
$F$.

(3.1.3) Assume that for any two general points   $y_1,y_2\in Y_{\bar
F}$  there is a morphism $f=f(y_1,y_2):\p^1\to  Y_{\bar F}$ such
that, $y_1,y_2\in
\im f$, 
$Y_{\bar F}$ is smooth along $\im f$ and $f^*T_{Y_{\bar F}}$ is semi
positive. Then any birational selfmap of $Y$ preserves $f$. Hence
there is an exact sequence
$$ 1\to \bir(Y_F)\to \bir(Y)\to \bir (X).
$$
\endproclaim

\demo{3.1.4 Remark}  The assumption (3.1.3) is satisfied if $Y_{\bar
F}$ is separably rationally connected (cf.
\cite{Koll\'ar96,IV.3.2}). More generaly, it also holds for the
cyclic covers of $\p^n$ that we are considering, \cite{ibid,V.5.19}. 
\enddemo

\demop   $M$ is big, hence there is an open set $U\subset X$ such
that sections of $M^k$ separate points of $U$ for $k\gg 1$. In
particular, if $g:C\to X$ is a nonconstant  morphism from a smooth
proper curve to $X$ whose image intersects
$U$, then
$\deg g^*M>0$.  

Let $g:C\to Y$ be a morphism such  that $g^*\Omega^1_Y$ is semi
negative. We have a map
$$ g^*h: g^*f^*M\to \wedge^i g^*\Omega^1_Y.
$$ Thus either $(f\circ g)(C)\subset X-U$ or  $(f\circ g)(C)$ is a
single point. This will allow us to identify the fibers of $f$.

In order to prove (3.1.1) pick a general point $z\in Z$ and let
$\phi_z:\p^1\to Y$ be the restriction of $\phi$ to $\{z\}\times
\p^1$.  
$$
\Omega_{Z\times \p^1}^1\vert \{z\}\times \p^1\cong
\o_{\p^1}^{n-1}+\o_{\p^1}(-2),
$$ and $\phi$ gives a map
$$
\Phi: (f\circ \phi_z)^*M@>{\phi_z^*h}>>
\phi_z^*\wedge^i\Omega_Y^1@>{\wedge^id\phi}>>
\wedge^i(\o_{\p^1}^{n-1}+\o_{\p^1}(-2)),
$$ which is nonzero for general $z$ since $\phi$ is separable. Thus
$\deg (f\circ \phi_z)^*M\leq 0$.  By the above remarks, this implies
that
$f\circ \phi_z$ is a constant morphism.

Pick a point $0\in \p^1$ and define $\psi:Z\to X$ by
$\psi(z)\deq f\circ \phi(z,0)$. This shows (3.1.1). 

Let $Z_F$ be the generic fiber of $\psi$. We obtain a dominant
$F$-morphism
$\phi_F:Z_F\times \p^1\to Y_F$ which is birational (resp. separable)
iff
$\phi$ is birational (resp. separable). 

Conversely, if $W_F$ is any variety and $W_F\to Y_F$  a morphism then
it extends to a map $W\map Y$ of the same degree. This shows (3.1.2).

Finally assume (3.1.3). Then there is an open set $Y^0\subset Y$ such
that if
$y_1,y_2\in Y^0$  and $f(y_1)=f(y_2)$  then  there is a morphism
$f=f_{(y_1,y_2)}:\p^1\to  Y$ such that, $y_1,y_2\in
\im f$, 
$Y$ is smooth along $\im f$ and $f^*T_Y$ is semi positive.

Let $\phi:Y\map Y$ be a  birational selfmap of $Y$; $\phi$ is defined
outside  a codimenion 2 set $Z\subset Y$.  By \cite{Koll\'ar96,
II.3.7}, the image of the general $f_{(y_1,y_2)}$  is disjoint from
$Z$. Thus we have an injection 
$$ f_{(y_1,y_2)}^*T_Y\hookrightarrow (\phi\circ f_{(y_1,y_2)})^*T_Y,
$$ which shows that the latter is also semi positive. Thus  
$\phi(y_1)$ and $\phi(y_2)$ are in the same fiber of $f$. Therefore 
$\phi$  preserves $f$, which gives the exact sequence 
$$ 1\to \bir(Y_F)\to \bir(Y)\to \bir (X).\qed
$$ 
\enddemo

Applying (3.1) to the projection $Y\to \p^m\times\p^n\to \p^m$ of
(2.5)
 we obtain: 

\proclaim{3.2 Corollary}    Fix a prime $p$ and let $k$ be an
algebraically closed field
 of characteristic 
$p$. Let $s\in H^0(\p^m\times\p^n,\o(a,b)^{\otimes p})$ be
 a general section. Assume that $ap-m-1>0$ and $bp-n-1=0$. 

Then 
$Y'\deq \p^m\times\p^n[\root{p}\of{s}]$ is ruled  iff  the generic
fiber of $Y'\to \p^m$ is  ruled  over the field $k(x_1,\dots,x_m)$.

Furthermore, 
$Y'$ has a degree $d$  separable uniruling  iff  the generic fiber of
$Y'\to \p^m$ has a degree $d$  separable uniruling over the field
$k(x_1,\dots,x_m)$.\qed
\endproclaim

(3.2) naturally leads to the following:

\demo{3.3 Question}  Let $X_F$ be a variety over a field. When is
$X_F$ ruled over $F$?
\enddemo

The problem is mainly interesting when $X_F$ is geometrically ruled.
I can not say much in general, so I consider only two simple
examples: conics and Del Pezzo surfaces. One should   keep in mind
that in our applications $F$ is the function field of a variety in
positive characteristic, so $F$ is not perfect. Also, we need these
results in  characteristic 2 for conics and in characteristic  3 for
cubic surfaces. These are the most unusual cases.

\demo{3.3.1 Example}  If $X_F$ is an arbitrary  variety which has no
$F$-points, then $X_F$ is not rational, but it can happen   that it
is ruled. For instance, if $Y_F$ has no $F$-points then $Y\times
\p^1$ is ruled and has no $F$-points.

This can happen even for quadrics in $\p^3$.  For example, choose
$a,b\in F$ such that $C=(x_0^2+ax_1^2+bx_2^2=0)$ has no
$F$-points (say $F=\r$ and $a=b=1$).  Then $C\times \p^1$ is
birational to the quadric 
$Q=(y_0^2+ay_1^2+by_2^2+aby_3^2=0)$ via the map
$$
\phi: (x_0:x_1:x_2,s:t)\mapsto (sx_0+atx_1:sx_1-tx_0:sx_2:tx_2).
$$
$Q$ is has no $F$-points.  \enddemo

\head 4.  Conics over nonclosed fields 
\endhead

The aim of this section is to study conics over arbitrary fields. We
study when they are ruled or uniruled. The main result is (4.1), but
for the applications we need (4.2).

\proclaim{4.1 Theorem}  Let $F$ be a field and $C_F\subset \p^2_F$ an
irreducible and reduced conic. The following are equivalent:

(4.1.1)  $C_F$ has a point in $F$.

(4.1.2) $C_F$ is ruled.

(4.1.3) $C_F$ has an odd degree uniruling.

\noindent If $C_F$ is smooth then these are also equivalent to

(4.1.4) $C_F\cong \p^1_F$. 
\endproclaim

\demop   Let $P$ be an $F$-point of  $C$. If
$C$ is geometrically irreducible, then projecting it from $P$ gives a
birational map
$C\to \p^1_F$, hence $C$ is ruled. Otherwise, projection exhibits $C$
as  a cone over a length two subscheme of $\p^1_F$, thus again  $C$
is ruled. This shows that  (4.1.1) $\Rightarrow$  (4.1.2) while 
(4.1.2) $\Rightarrow$  (4.1.3) is clear.

The proof of (4.1.3) $\Rightarrow$  (4.1.1) is longer.  Let $A$ be a
zero dimensional $F$-scheme and $\phi_F:A\times_F\p^1\to C_F$ an odd
degree uniruling. Let $A_i\subset A$ be the irreducible components.
$\deg
\phi_F=\sum_i\deg (\phi_F|A_i\times \p^1)$, thus one of the $\deg
(\phi_F|A_i\times
\p^1)$ is odd. Thus we may assume that $A$ is irreducible.  
$\deg (\phi_F|\red A\times \p^1)$ divides $\deg \phi_F$, hence we may
also assume that $A=\spec_FF'$ where $F'\supset F$ is a  field
extension.

By base change to $F'$  we obtain
$$
\p^1_{F'}\hookrightarrow \spec_{F'}(  F'\otimes_FF')\times \p^1_F\to
C_{F'}.
$$ Thus by the L\"uroth theorem, every irreducible component of $\red
C_{F'}$ is birational to $\p^1_{F'}$.  Passing to the algebraic
closure we obtain
$$
\phi_{\bar F}:\spec_{\bar F}(\bar F\otimes_FF')\times \p^1\to C_{\bar
F}.
$$ The left hand side may have several irreducible components,  
conjugate to each other. Let $\bar \phi:\p^1\to \p^1$ be the induced
map between  any of the irreducible and reduced components. By
counting degrees we obtain the following:

\proclaim{4.1.5 Claim} Notation as above. Then  

(4.1.5.1) $\deg \phi_F=\deg(F'/F)\deg (\bar \phi)$ if $C_{\bar F}$ is
a smooth conic, and

(4.1.5.2) $2\deg \phi_F=\deg(F'/F)\deg (\bar \phi)$ if $C_{\bar F}$
is a singular conic.

Thus $\deg (F'/F)$ is odd in the first case and is not divisible by 4
in the second case.\qed
\endproclaim

Assume first that $C_{\bar F}$ is a smooth conic. 
$A\times_F\spec F'$ has a closed point, thus we get a morphism
$\p^1_{F'}\to C_{F'}$. Hence $C_{F'}$ has a point in $F'$. This in
turn gives an odd degree point on $C_F$; let 
 $L$ be the corresponding line bundle. The restriction of
$\o_{\p^2}(1)$ to
$C_F$ is a line bundle of degree 2. We conclude that $C_F$ has a line
bundle of degree 1. Any of its sections gives an $F$-point  on $C_F$.

If $C_{\bar F}$ is a pair of intersecting lines, then the
intersection point is defined over $F$.

Finally consider the case when $C_{\bar F}$ is a double line; this can
happen only in  characteristic 2.   The equation of $C_F$ is $\sum
b_ix_i^2=0$.

\proclaim{4.1.6 Claim} Let $C_F=(\sum b_ix_i^2=0)$ be an irreducible
conic over a field of characteristic 2.  Let $E/F$ be a separable
extension. Then any $E$-point of
$C_F$ is an
$F$-point.
\endproclaim

\demop  We may assume that $E/F$ is Galois. Assume that    $P$ is an
$E$-point
 which is not an $F$-point.  Conjugates of
$P$ over
$F$ also give $E$-points, thus we obtain that $\red C_{\bar F}$ is
defined over
$E$. $\red C_{\bar F}$ is also defined over the purely inseparable
extension 
$F^i=F(\sqrt{b_0},\sqrt{b_1},\sqrt{b_2})$, hence also over the
intersection
$F''\cap F^i=F$ (cf. \cite{Koll\'ar96, I.3.5}). This is a
contradiction. 
\qed\enddemo

As in the irreducible case we know that $C_{F'}$ has an $F'$-point.
Thus by (4.1.6) 
$F'/F$ is not separable. In view of (4.1.5.2) we conclude that  there
is a subextension $F'\supset F''\supset F$ such that $F'/F$ has odd
degree (hence separable) and $F'=F''(\sqrt{s})$ for some $s\in F''$. 
By (4.1.6) it is enough to show that $C$ has an $F''$-point. 

As we mentioned earlier, $\red C_{F'}$ is birational to $\p^1_{F'}$,
thus 
$\red C_{F'}$ is a line in $\p^2$. Therefore 
 $C_{F'}$ is a double line with equation
$(\sum a_ix_i)^2=0$ where $a_i\in F'$.   $a_i^2\in F''$, thus the
equation of $C$ over $F''$  is $\sum a_i^2x_i^2=0$.  Write
$a_i=c_i+sd_i$ where $c_i,d_i\in F''$. The equation of $C$ is
$$
\sum a_i^2x_i^2=(\sum c_ix_i)^2+s^2(\sum d_ix_i)^2=0.
$$ The  solution of $\sum c_ix_i=\sum d_ix_i=0$ gives an $F''$-point
on
$C$.

The equivalence with (4.1.4) was established in the course of the
proof. \qed\enddemo

\demo{4.1.7 Remark}  If $C_{\bar F}$ is a smooth conic, then
$C_F$ has a degree 2  separable uniruling. Take any general line in
$\p^2$. Its intersection points with $C_F$ are in a degree 2
separable extension of
$E\supset F$. Thus $C_E\cong \p^1_E$.
\enddemo

It remains to establish that the generic fibers apearing in (3.1) are
not ruled. There should be a general result about conics over function
fields, but I could not find a simple   proof. For our applications
the following is sufficient:

\proclaim{4.2 Proposition}   Let $k$ be a  field  of characteristic 2
and
$F=k(x_1,\dots,x_m)$  the field of rational functions in $m\geq 2$
variables. Fix an even integer $d\geq 2$ and let $a,b,c\in
k[x_1,\dots,x_m]$ be   general polynomials of degree $d$. Then the
conic
$$ C=(y_0^2=ay_1^2+by_1y_2+cy_2^2)\subset
\p^2_F
$$ is not ruled (over   $F$). Moreover $C$ does not have any odd
degree uniruling.
\endproclaim

\demo{4.2.1 Comments}  It is worth while to remark that for (4.2) to
hold it is essential  that  $k(x_1,\dots,x_m)$ is not perfect. A
point on $C$ is given by 
$P=(\sqrt{a}, 1,0)$. If $F$ is perfect of characteristic 2, then $P$
is an $F$-point and $C$ is rational.
\enddemo

\demop By (4.1) it is sufficient to establish that $C$ has no
$F$-points.  We can identify $ay_1^2+by_1y_2+cy_2^2$ with a section
$s$ of
$\o_{\p^{m-1}\times \p^1}(d,2)$. Let $Y\deq \p^{m-1}\times
\p^1[\sqrt{s}]$ be the corresponding double cover. The generic fiber
of $\pi: Y\to \p^{m-1}$ is
$C$, thus it is sufficient to prove that $\pi$ has no sections.  By
(2.2.3) $Y$ has only isolated singularities. The fibers  of $\pi$ over
$b=0$ are double lines. This shows that $\pi$ does not even have local
sections at the generic point of $b=0$.\qed\enddemo

\head 5. Del Pezzo surfaces over nonclosed field
\endhead

In order to complete the proof of (1.3) we need to study the Del
Pezzo surfaces:

(5.1.1) $S_F$ with equation $u^3=f_3(x,y,z)$ in characteristic  3.

\noindent Although we do not need it, it is very natural to study
also the surfaces:

(5.1.2)  $T_F$  with equation $u^2=f_4(x,y,z)$ in characteristic  2.

\noindent  To get an idea of the geometry of these surfaces, we study
them first over perfect fields.

\demo{5.2 Remark: The case of perfect fields}  Assume first that our
base field is algebraically closed. It is easy to see that  we can
write
$$ f_3=l_1l_2l_3+l_4^3,\qtq{and} f_4=l_1l_2l_3l_4+q^2
$$ where the $l_i$ are linear and $q$ quadratic. Thus we can make
coordinate changes $u=u-l_4$ (resp. $u=u-q$) and then a suitable
coordinate change among the $x,y,z$  to reduce the equations to the
form
$$ u^3=xyz,\qtq{and} u^2=xyz(x+y+z).
$$ From this we see that the cubic $S$ has three singular points of
type
$A_2$.  The degree two Del Pezzo $T$ has seven singular points of
type $A_1$  at the seven points of  $\Bbb F_2\p^2$.

Next consider the case when the base field $F$ is perfect. Then
these   singular points are defined over
$F$.  For the cubic $S_F$   resolve the singular points, and contract
the birational transforms of the 3 coordinate axes to get a degree 6
Del Pezzo surface. Over a perfect field $F$ a  degree 6 Del Pezzo
surface is rational iff it has an $F$-point
\cite{Manin72,IV.7.8}.   Our surface $S_F$ does have  $F$-points over
perfect fields of characteristic 3, for example
$P=(1,0,0,\root{3}\of{f_3(1,0,0)})$.

 For the surface $T_F$   resolve the singular points,  and contract
the birational transforms of the seven lines in $\Bbb F_2\p^2$ to get
a Brauer-Severi variety. It has a point in a  degree 7 extension,
hence it is isomorphic to $\p^2$.  Thus our surfaces $S_F$ and $T_F$
are rational over any perfect field. 
\enddemo

In our case the field $F$ is not perfect, the surfaces $S_F$ and
$T_F$ are nonsingular over $F$ but they are not smooth. In order to
show that they are not ruled, we have to understand how the presence
of nonsingular but nonsmooth points effects the birational geometry
of a surface.

The crucial result is the following:

\proclaim{5.3 Theorem}  Let $F$ be a field and $S,T$ nonsingular and
proper surfaces over $F$. Assume that $T$ is smooth except possibly
at finitely many points. Let $f:S\map T$ be a birational map.

Then $f$ is defined at all nonsmooth points of $S$.
\endproclaim

\demop  There is a sequence of blowing ups of closed points $p:S'\to
S$ such that $f\circ p:S'\to T$ is a morphism.  Let $P\in S$ be a
closed nonsmooth point, and assume that
$f$ is not defined at $P$.  Then $p^{-1}(P)$ is 1-dimensional and
there is an irreducible component
$E\subset p^{-1}(P)$ such that $p\circ f$ is a local isomorphism at
the generic point of $E$. This implies that $S'$ is smooth at the
generic point of
$E$. This contradicts the following lemma.\qed\enddemo

\proclaim{5.3.1 Lemma}  Let $F$ be a field and $S, S'$ nonsingular  
surfaces over $F$. Let $p:S'\to S$ be a proper and birational
morphism and
  $P\in S$  a closed nonsmooth point.

 Then $S$ is not smooth at all points of $p^{-1}(P)$.
\endproclaim

\demop By induction it is sufficient to consider the case when $p$ is
the blow up of $P$. We may also assume that  $S$ is an affine
neighborhood of $P$ such that the maximal ideal $m_P$ is generated by
two global sections $u,v\in m_P\subset \o_S$. The blow up can be
described by two affine charts, one of them is
$$ S'_1\deq (u-vs=0)\subset S\times \a^1, \qtq{where $s$ is a global
coordinate on
$\a^1$.}
$$ By assumption $S$ is not smooth at $P$, so $S\times \a^1$  is not
smooth along $P\times \a^1$.  $S'_1$ is a Cartier
 divisor on $S\times \a^1$, thus is it also not smooth along $P\times
\a^1$. This was to be proved.\qed\enddemo

In the course of the proof we used the following  results about
birational transformations of nonsingular surfaces. For a proof see,
for instance, \cite{Zariski58}.

\proclaim{5.3.2 Proposition}   Let $S,T$ be proper and nonsingular
surfaces over a field $F$ and $\phi:S\map T$ a birational map. Then
$\phi$ is a composite of blow ups and blow downs (of closed points). 
In particular,
$h^i(S,\o_S)=h^i(T,\o_T)$. \qed
\endproclaim

\proclaim{5.3.3 Corollary}  Let $F$ be a field and $S$ a
nonsingular   and proper surface over $F$.  Assume that $S$ is smooth
except possibly at finitely many points. Let $f:S\map S$ be a 
birational map. 

 Then $f$ is a local isomorphism at all  nonsmooth  points of
$S$.
\endproclaim

\demop Factor $f$ as
$$ f:S@<p<< S'@>p'>> S
$$ where $p,p'$ are birational morphisms.  By (5.3)  we see that 
 $p$ and $p'$ are both local isomorphisms at the nonsmooth  points of
$S$.\qed\enddemo

This gives the following rationality criterion:

\proclaim{5.4 Corollary}  Let $F$ be a  field and $S$ a nonsingular  
and proper surface over $F$. Assume that $S$ is generically smooth.
The following are equivalent:

(5.4.1) $S$ is rational (over $F$).

(5.4.2) There is a  two dimensional linear system $L=|C|$ on $S$ with
(infinitely near) base points $P_i$ of multiplicity $m_i$ such that

(5.4.2.1) a general $C\in L$ is birational to $\p^1$;

(5.4.2.2) $S$ is smooth along a general $C\in L$;

(5.4.2.3) $C\cdot K_S+\sum m_i=-3$ and $C^2-\sum m_i^2=1$.
\endproclaim

\demop Assume that there is a birational map $f:S\map \p^2_F$. Let
$Z\subset S$ be the locus of nonsmooth points. By (5.3), $f$ is
defined along $Z$ and
$f(Z)$ is zero dimensional. Set
$L=f^{-1}_*|\o_{\p^2}(1)|$.  (5.4.2.1--2) are clear,  and (5.4.2.3)
is the usual equalities (5.4.3.1).

Conversely, assume (5.4.2). Resolve the base points of $L$ to obtain
a base point free linear system $L'$ on $S'$. 
 From (5.4.2.3) we obtain that
$$ C'\cdot K_{S'}=-3\qtq{and}{C'}^2=1\qtq{ for $C'\in L'$.}
$$ Thus  the linear system $L'$  maps $S'$ birationally to $\p^2_F$. 
\qed\enddemo

What we really need is a ruledness criterion. If $S$ is ruled, it can
be birational to a surface which is the product of $\p^1$ with a
conic which has no $F$-points. Thus the natural linear system
obtained on $S$ is two dimensional and its general member is
geometrically reducible. I found it   clearer to concentrate instead
on a single curve, which may not be defined over
$F$. We do not get an equaivalence any longer, but for our
applications this does not matter. 

\demo{5.4.3 Remark} For linear systems we considered base points with
multiplicities. If we look at a general member, the corresponding
notion is a curve with assigned multiplicities. All infinitely near
multiple points are assigned with their multiplicity, but we may also
have some smooth points assigned with multiplicity one. 

If $C\subset S$ is a curve on a smooth surface with assigned
multiplicities
$m_i$ at the points
$P_i$ and $\phi:S\to S'$ is a birational map, there is a natural
birational transform  $C'\subset S'$  with assigned multiplicities
$m'_i$ at the points
$P'_i$. In order to define this we need only the case when $\phi$ is
the blow up of a point or its inverse, where the definition is the
obvious one. The values of the expressions
$$ C\cdot K_S+\sum m_i\qtq{and}C^2-\sum m_i^2\tag 5.4.3.1
$$ are birational invariants, cf. \cite{Hudson27, p.5}.
\enddemo

\proclaim{5.5 Corollary}  Let $F$ be a field and $S$ a nonsingular  
and proper surface over $F$ such that $S$ is generically smooth,
geometrically irreducible  and   $h^1(S,\o_S)=0$.  Assume that $S$ is
ruled (over $F$).

Then there is a rational curve $C\subset S_{\bar F}$ with assigned
(infinitely near) multiple points $P_i$ of multiplicity $m_i$ such
that

(5.5.1) $S_{\bar F}$ is smooth along $C$.

(5.5.2) $C\cdot K_S+\sum m_i=-2$ and $C^2-\sum m_i^2=0$.

(5.5.3) $\o_{S_{\bar F}}(2C)\in \pic (S)$. 
\endproclaim

\demop By asumption there is a nonsingular, geometrically integral
 curve $D$  and  a birational map $f:S\map D\times \p^1$.  By
(5.3.2), 
$$ 0=h^1(S,\o_S)=h^1(D\times \p^1,\o_{D\times \p^1})=h^1(D,\o_D).
$$ Thus $D$ is isomorphic to a smooth conic. 

Let $d\in D_{\bar F}$ be a general point, $C'= d\times
\p^1$ and  $C=f^{-1}_*(C')\subset S_{\bar F}$  the corresponding
birational transform.  $S$ is smooth along  $C$ by (5.3), and (5.5.2)
is the usual equalities (5.4.3.1).

Finally, $D$ has a degree 2 point defined over $F$, thus the line
bundle
$\o_{D\times \p^1}(2C')$ is defined over $F$. Therefore $\o_{S_{\bar
F}}(2C)$ is also defined over $F$. \qed\enddemo

The main result of this section is (5.7). Over perfect fields it is a
special case of more general results of Segre and Manin (see
\cite{Manin72}).  The proofs in \cite{Manin72} use  the structure of
the Picard group of smooth Del Pezzo surfaces  to compute the action
of certain involutions. In our case the Picard groups are small (5.6),
and it is easier to use the geometric ideas of 
\cite{Segre43,51} to analyze the involutions. (5.3) essentially says
that all the relevant geometry takes place inside the smooth locus,
where the geometric description of the involutions works well.

\proclaim{5.6 Lemma} Let $S_F$ denote an integral cubic  surface
$u^3=f_3(x,y,z)\subset \p^3$ for $\chr F=3$, or an integral double
plane
 with equation $u^2=f_4(x,y,z)$ for $\chr F=2$. 

Then $\pic(S_F)=\z[-K_S]$.
\endproclaim

\demop Let $\pi:S_F\to \p^2$ be the projection to the $(x,y,z)$-plane.
$\pi$ is purely inseparable, thus if $C\subset S_F$  is any divisor
then
$\pi^*(\pi_*(C))=(\deg \pi)C$. Therefore $-K_S=\pi^*(\o_{\p^2}(1))$
generates $\pic(S_F)\otimes \q$. Since $(K_S^2)\leq 3$, we get that
$K_S$ is not divisible in $\pic(S_F)$, hence
$-K_S$ generates $\pic(S_F)/(\text{torsion})$. Thus we are left to
prove that
$\pic(S)$ has no  torsion.

Let $[C]\in \pic(S)$ be a numerically trivial Cartier divisor. By
Riemann-Roch,
$$ h^0(\o_S(C))+h^0(\o_S(K_S-C))\geq \chi(\o_S)=1.
$$
$-K_S$ is ample, so $h^0(\o_S(K_S-C))=0$. Thus $\o_S(C)$ has a
section and
$\o_S(C)\cong \o_S$.\qed\enddemo

\proclaim{5.7 Theorem}  Let $F$ be a field and $S$ a  nonsingular  
Del Pezzo surface over $F$. Assume that $S_{\bar F}$ is integral, 
$\chi(\o_S)=1$, 
$\pic(S)=\z[-K_S]$ and $K_S^2=1,2,3$.

Then $S$ is not ruled  (over
$F$).  
\endproclaim

\demo{5.7.1 Remark} Let $S$ be a  nonsingular   Del Pezzo surface
over $F$ such that $S_{\bar F}$ is integral, 
$\chi(\o_S)=1$, 
 and $K_S^2=1,2,3$. The  structure theory of Del Pezzo surfaces 
 shows that $S$ is a cubic surface for $K_S^2=3$, a double plane for
$K_S^2=2$ and as expected for $K_S^2=1$ (cf. \cite{Koll\'ar96,
III.3}). 

\cite{Reid94} contains examples of nonnormal Del Pezzo surfaces  with
$\chi(\o_S)<1$.  Some of these may have nonsingular models over
nonperfect fields. I have not checked if this indeed happens or what
can be said about their arithmetic properties. 
\enddemo

\demop  Assuming that $S$ is ruled, we derive a contradiction. (5.5)
guarantees the existence of a rational curve $C\subset S_{\bar F}$
satisfying the properties (5.5.1--3).  By (5.5.3) $\o(2C)$ is in
$\pic(S)$. Since $K_S^2\leq 3$, $K_S$ is not divisible by 2 in
$\pic(S_{\bar F})$ and therefore
$\o(C)$  is in $\pic(S)$. This implies that 
$C\in |-dK_{S_{\bar F}}|$ for some $d$.  We   show that there can not
be such a curve:

\proclaim{5.8 Proposition} Let $S$ be an integral Del Pezzo surface
over an algebraically closed field such that $\chi(\o_S)=1$, 
 and $K_S^2=1,2,3$.  Then
$S$ does not contain any curve $C$ satisfying the following
properties:

(5.8.1) $C$ is birational to $\p^1$;

(5.8.2) $S$ is smooth along $C$;

(5.8.3)  $C\in |-dK_S|$ for some $d$.

(5.8.4)  $\sum_i m_i=d(K_S^2)-2$ and $\sum_i m_i^2=d^2(K_S^2)$, where
 $P_i$ are the (assigned  and infinitely near) multiple points of $C$
with multiplicity
$m_i$.
\endproclaim

\demop If $m_i\leq d$ for every $i$ then
$$
\sum_i m_i^2\leq d\sum_i m_i=d^2(K_S^2)-2d<d^2(K_S^2),
$$
 a contradiction. Thus there is a point $P=P_1$ such that $m=m_1>d$.

Let $P\in L\subset S$ be a line (that is, $-K_S\cdot L=1$). Then
 $m\leq (C\cdot L)=d$.  Therefore $P$ can not lie on any line.

If $(K_S^2)=1$ then any member of $|-K_S|$ is a line, thus there is a
line through any point. Hence we are done if $(K_S^2)=1$.

Next consider the case when $(K_S^2)=3$. That is, $S\subset \p^3$ is
a cubic surface.

Let
$D\subset S$ be the intersection of $S$ with the tangent plane at
$P$.  Since there is no line through $P$, $D$ is an irreducible cubic
whose unique singular point is at $P$. In particular, $D$ is
contained in the smooth locus of $S$. 

The point $P$ determines a birational selfmap $\tau$ of $S$ as
follows. Take any point $Q\in S$, connect $P,Q$ with a line and let
$\tau(Q)$ be the third intersection point of the line with $S$.
$\tau$ is an automorphism of $S-D$. Another way of describing $\tau$
is the following. Projecting $S$ from $P$ gives a diagram
$$
\CD B_PS@>q>> \p^2\\ @VpVV @.\\ S @.
\endCD
$$ Let $E\subset B_PS$ be the exceptional curve; $C', D'\subset B_PS$
the birational transforms of $C,D$. $q$ is a degree two morphism and
$\tau$ is the involution interchanging the two sheets. (Computing
with the local equation at $P$ shows that $q$ is separable  in
characteristic 2 if $D$ is irreducible, thus $\tau$ always exists.) 
Furthermore, $\tau(E)=D'$ and 
$p^*\o_S(1)(-E)=q^*\o_{\p^2}(1)$. Thus
$$ C'+(m-d)E\in |q^*\o_{\p^2}(d)|,\qtq{hence} 
\tau(C'+(m-d)E)\in |q^*\o_{\p^2}(d)|.
$$ Pushing this down to $S$ we obtain that 
$$
\tau(C)+(m-d)D\in |\o_S(d)|,\qtq{hence}  \tau(C)\in |\o_S(d-(m-d))|.
$$ Thus $\tau(C)$ satisfies all the properties (5.8.1--4) and its
degree is lower than the degree of $C$. We obtain a contradiction by
induction on
$d$.

A similar argument works if $(K_S^2)=2$, but the details are a little
more complicated. I just outline the arguments, leaving out some
simple details. 

We already proved that $P$ is not on any line, and a similar argument
shows that $P$ can not be a singular point of a member of $|-K_S|$.

Since $h^0(S,-2K_S)=7$, there is a curve $D\in |-2K_S|$ which has a
triple point at $P$. In fact, $D$ is unique, it is a rational curve
and $P$ is its only singular point. Thus $S$ is smooth along $D$.  As
before we look at the blow up diagram
$$
\CD B_PS@>q>> Q\subset \p^3\\ @VpVV @.\\ S @.
\endCD
$$ where $Q$ is a quadric cone; the image of $B_PS$ by the linear
system
$|-2K_{B_PS}|$. $q$ is a degree two morphism and $\tau$ is the
involution interchanging the two sheets. (Again one can see  that $q$
is separable  in characteristic 2 if $|-K_S|$ does not have a  member
which is singular at $P$.) Let $E\subset B_PS$ be the exceptional
curve; $C', D'\subset B_PS$ the birational transforms of $C,D$.   
 $\tau(E)=D'$ and  $p^*\o_S(2)(-2E)=q^*\o_{\p^3}(1)$. As before we
obtain that 
$$
 \tau(C)\in |\o_S(d-2(m-d))|.
$$
 We obtain a contradiction by induction on
$d$.
\qed\enddemo

\vskip1cm

University of Utah, Salt Lake City UT 84112 

kollar\@{}math.utah.edu

\end